\documentclass[12pt,twoside, reqno]{amsart}
\usepackage[latin1]{inputenc} 
 \usepackage[english]{babel}
\usepackage{amsthm}
\usepackage{amssymb, verbatim}
\usepackage{amsmath} 
\usepackage{eucal}
\usepackage{mathrsfs} 
\usepackage{hyperref}

 \topmargin = - 30 pt %distanza dal margine superiore 0= un bel po` negativo diminuisce
\usepackage[all]{xy}
\usepackage{amsmath,amsthm,amssymb,amsopn,amscd,tikz-cd}
\newtheorem{thm}{Theorem}[section]

\theoremstyle{definition}

\newtheorem{defn}[thm]{Definition}

\numberwithin{equation}{section}

\begin{document}

\title[The braided Doplicher-Roberts program]{
 The braided Doplicher-Roberts program and the Finkelberg-Kazhdan-Lusztig equivalence: A historical perspective,  recent progress, and 
future directions}

 \author[C.~Pinzari]{Claudia Pinzari }
 \email{ 
 pinzari@mat.uniroma1.it}
\address{Dipartimento di Matematica, Sapienza Universit\`a  di
Roma\\ P.le Aldo Moro, 5 -- 00185 Rome, Italy}\maketitle

\centerline\small{\it Dedicated to the memory of Sergio Doplicher}

 \begin{abstract}   
 Our recent approach to the Finkelberg-Kazhdan-Lusztig equivalence theorem centers on the construction of a fibre functor associated with the categories in the equivalence theorem,
  which in turn explains the underlying algebraic and analytic structure of  the corresponding weak Hopf algebra
in a new sense. 
We provide a non-technical and historical overview of the core arguments behind our   proof, discuss these structural properties,  and its applications to rigidity and unitarizability   of braided fusion categories arising from conformal field theory. We conclude proposing  some natural directions for future research.

  \end{abstract} 
  
    \tableofcontents

\section{Introduction}\label{Introduction}

Doplicher and Roberts originally proposed  extending the duality theory for compact groups of \cite{DR1},  \cite{DR_Annals} from rigid symmetric tensor $C^*$-categories with simple tensor unit
  to   unitarily  braided  categories \cite{DR_CMP}.   
  
 Their   approach in the symmetric case centers on constructing a unique faithful symmetric tensor functor to the  category of Hilbert spaces. By constructing  a unique fibre functor from a minimal rigid symmetric tensor subcategory generated by a single object, and then extending it to the full subcategory generated by the object via a crossed product construction, they  showed commutativity
 of a certain relative commutant $C^*$-algebra in the crossed product. In this way they
  constructed a classical action of the compact Lie group ${\rm SU}(d)$ obtained by Schur-Weyl duality, by homeomorphisms on the Gelfand spectrum of the commutative algebra, with $d$ the statistical dimension of the generating object, automatically a positive integer  as recovered from rigidity and permutation symmetry.
By applying a Mackey induction from this commutative $C^*$-algebra  they construct a compact Lie group as a closed subgroup such that its representation category is equivalent to the given category, and  unique up to conjugation. If the category is not singly generated, they obtain a   compact group by an inductive limit procedure unique up to isomorphism. The symmetric nature of the functor is the key property ensuring the group's existence and uniqueness.

The main application of their theory was given in the axiomatic setting of algebraic quantum field theory (AQFT), founded by
Haag and Kastler
\cite{HK64}, 
where Doplicher and Roberts were able to construct
  a compact gauge group of internal symmetries, the physical Hilbert space, and the field algebra
   directly from the local observable algebras and their localized endomorphisms \cite{Haag}.
   This approach ensured that   quantum fields automatically satisfy   Bose-Fermi statistics and that the conjugate 
   localized endomorphisms, representing particle-antiparticle duality, automatically exist \cite{DR_CMP}.
  
  Extending the Doplicher-Roberts reconstruction theorem to braided categories is interesting for the theory of quantum fields. While
 rigorous    models      in low-dimensional QFT such as conformal field theory (CFT),
are well-developed, a general framework to reconstruct the field algebra and a quantum gauge group is still emerging.
This program, which requires new methods to overcome the obstacles posed by braided symmetry, can be described as {\it the braided Doplicher-Roberts program.}

      In this paper we provide an overview of   some of our recent work  \cite{FKL} regarding the braided Doplicher-Roberts
    program,   by situating  it within a historical context and propose natural directions for future development.

 In the case of braided categories,   it soon became clear that  the statistics allows more than Bose-Fermi alternative.
 Indeed, in algebraic CFT, the statistical dimension is often non-integer, see \cite{Haag}.
    
 Parallel to these developments, the search for a self-consistency program   to solve quantum field theories through their internal symmetries
 --pioneered by Polyakov \cite{Pol70}--found its    $2D$ realization in
  what is now known as the Wess-Zumino-Witten (WZW) model. The current algebra of this model was shown by Witten \cite{Witten1984}
     to satisfy the commutation relations of affine Kac-Moody algebras at positive integer levels  \cite{Kac2}. In the same paper, Witten   showed the equivalence between the WZW model and a theory of free fermions in 2D. At the same time, Belavin, Polyakov and Zamolodchikov (BPZ) established an axiomatic framework
     of 2D CFT by introducing the   Virasoro algebra and {\it primary fields} \cite{BPZ}.
  Finally,   Knizhnik and Zamolodchikov connected the current algebra of Witten to the work by BPZ
  by showing that    the $n$-point correlation functions of primary fields of the WZW model  satisfy
  a system of partial differential equations, the KZ equations   \cite{KZ84}.

  Frenkel and Zhu established a vertex operator algebra structure on a certain representation $L_{k, 0}$ of the
  affine Kac-Moody (with $k$ the level)   and analogously for the Virasoro algebra  \cite{Frenkel_Zhu}, in such a way that the representations of $L_{k, 0}$ give rise to the WZW model.
  The complete operator algebraic construction  of a conformal net associated to the corresponding
   loop group--the net of observable algebras in the setting of chiral conformal AQFT--was later provided by Wassermann who
 described the fusion of   representations  in terms of bimodules over von Neumann algebras and Connes' relative
tensor product structure
 \cite{Wassermann}. Finally, the systematic braided tensor equivalence between the representation category of the unitary affine vertex operator algebras (endowed with the Huang-Lepowsky tensor structure)  and the DHR categories of the corresponding loop group conformal nets was established by Gui in a series of papers  \cite{GuiII}, \cite{GuiIII}, \cite{GuiV}.

While these $2D$ models are now well-understood, Polyakov's original vision for $D\geq3$
has seen a renewed interest through modern numerical and analytical bootstrap methods \cite{Slava}. 
Whereas Polyakov's approach starts with local field operators and the operator product expansion (OPE), the Haag-Kastler framework seeks to characterize the theory through the structural properties of its local observables \cite{Haag}. 
%% The Doplicher-Roberts program   provides the rigorous reconstruction of the fields and
%% symmetries that Polyakov's bootstrap assumes to be present.

  %%The original  $4D$  Doplicher-Roberts theory is restricted to permutation symmetry and compact group duals, leaving the reconstruction of %%theories with non-integer statistical dimensions or braided symmetries--such as those encountered in 
%%  low-dimensional CFT 
%%  %%or in  the algebraic $4D$  QED model, which is described by string-localized charges \cite{MRS}
 %% --as an open problem.
%%The reconstruction of the field algebra and its underlying (quantum) gauge symmetry in these lower-dimensional contexts remains a %%fundamental challenge for the braided Doplicher-Roberts program.

       {\it Quantum groups} made their original appearance  in quantum integrable systems and models of low dimensional algebraic quantum field theory (conformal field theory), in connections with the discovery of
  a   method for studying the {\it quantum inverse scattering problem}. This method is based on the use of the quantum $R$-matrix and was first applied by Sklyanin in the late 70s,     inspired by the breakthrough works of Baxter in statistical mechanics of the early 70s.
          In the following section we give a brief historical background on quantum groups.

              Thus there arose the great interest of recognizing quantum group structure in conformal field theory.
           %% The monodromy of the KZ equations was shown by  Tsuchiya and Kanie  \cite{TK87}-------
            Tsuchiya and Kanie \cite{TK87} provided a rigorous treatment of the work by Knizhnik and Zamolodchikov. They considered the Kac-Moody algebra $\hat{\mathfrak sl_2}$ and proved
        that primary fields exist and are uniquely determined by their zero modes (certain intertwiners between representations
        of ${\mathfrak sl}_2$)  and characterized their values by   Clebsch-Gordan fusion rules for $\hat{\mathfrak sl}_2$.
            They showed that    the   representation of the braid group arising from the monodromy of the KZ equations   corresponds
             to representations of the Hecke algebras at a suitable root of unity as found by Jones and Wenzl  \cite{Wenzl_thesis}, \cite{Jones_braid},  \cite{Wenzl_Hecke}, originating from the theory of subfactors \cite{Jones_1}, and deeply related to the theory of knots  \cite{Jones_3}.
            Kohno       \cite{Kohno}  showed that for a non-exceptional complex simple Lie algebra ${\mathfrak g}$ and its vector representation,
            the monodromy representations of the braid group  are   equivalent to those arising from the $R$-matrix of  quantum deformations of Lie groups. Moreover he applied the results to CFT and obtained   unitary representations  from the analytic structure of subfactors or Hecke algebras.  
      
      Kohno's equivalence  was further studied  by Drinfeld for a formal deformation parameter $h$.  He constructed a   braided quasi-Hopf algebra structure on $U({\mathfrak g})$  and associated  to it a braided tensor category of representations--the Drinfeld category--with braiding and associativity morphisms derived from the KZ equations. Drinfeld showed that this category  is equivalent to the representation category of the formal quantization  $U_h({\mathfrak g})$  \cite{Drinfeld_quasi_hopf}.
     
        The full construction of a braided tensor category structure on the module category of affine Lie algebras at positive integer level
         and its
        equivalence with the fusion categories of quantum groups at roots of unity, 
        emerged from the combined work  of several authors in different settings, mainly  Kazhdan and Lusztig in an algebraic approach, who considered affine Lie algebras with non-positive rational levels,    
        Huang and Lepowsky who considered affine vertex operator algebras at positive integer levels, with an analytic approach based on differential equations, and  Finkelberg, with a geometric algebraic approach, based on work by Beilinson, Feigin and Mazur and Faltings, and Teleman. These works were originally influenced by the insight
        of Moore and Seiberg \cite{Moore_Seiberg2} in physics regarding the braided tensor category structures arising from rational CFT. 
        The part of the equivalence with the quantum group fusion categories is due to Kazhdan and Lusztig. The reader may refer to the introduction of \cite{Pinzari_constructing} for more historical information on this problem.
        
       As we shall discuss, the existing proof of this theorem
        is considered remarkably difficult and indirect. It is within this context that our work begins, seeking a unification
        through   the braided Doplicher-Roberts program.

        %%A rather recent proof of Drinfeld-Kohno theorem for a real deformation parameter has been given by 
        %%Neshveyev and Tuset in the setting of compact quantum groups and quasi-Hopf $C^*$-algebras. The equivalence in case of affine
         %%Lie algebras on positive integer levels on one side and fusion categories of quantum groups at roots of unity  is substantially more 
         %%complicated, has been studied by the joint works by Kazhdan and Lusztig, and Finkelberg, and we shall come back to this.  

                      Within the braided Doplicher-Roberts program, the equivalence amounts to asking whether the gauge symmetries arising from the CFT are described in a natural way by some kind of quantum group structure on a semi-simple algebra--or a suitable generalization endowed with suitable analytic and algebraic structure--naturally associated to the operator algebraic approach to CFT. Mack and Schomerus suggested   the search for an analogue of the global gauge group in the setting of Doplicher-Roberts program as a suitable truncation of a quantum group at roots of unity to eliminate non-physical representations arising  from  non-semisimplicity of $U_q({\mathfrak sl}_2)$
                      \cite{MS}.
                      
 A main result of  our work in \cite{FKL} is the construction of a {\it quantum gauge group} from the braided tensor category of modules of affine vertex operator algebras at positive integer level, naturally realized as a {\it unitary weak Hopf $C^*$-algebra} in a new sense, associated to the {\it zero modes of primary fields}. The structure of this quantum gauge group is summarized in Sect. 8 with the term of {\it unitary coboundary weak Hopf $C^*$-algebra}. We next discuss the genesis of our quantum group.

        Extending the notion of a symmetric functor to the braided setting remained a long-standing challenge, primarily because the only braided symmetry admitted by the category of vector spaces is the standard permutation symmetry. 
        
        Compact quantum groups in the operator algebraic setting were pioneered by Woronowicz \cite{Wor}, \cite{Wor2}, \cite{WoronowiczTK}, \cite{Wor_compact_quantum_groups}.

In our previous works, we extended Doplicher-Roberts duality to rigid tensor $C^*$-categories with simple tensor unit not necessarily braided.   More precisely, in  \cite{PRergodic} we    constructed a non-commutative $C^*$-algebra endowed with an ergodic action of a compact quantum group.  Our noncommutative  $C^*$-algebra was a quantum homogeneous space  playing the role of Doplicher-Roberts classical homogeneous space, and the ergodic action of the quantum group
  that of the induced representation of a closed  subgroup. This result met the property   discovered by Podles that for non commutative $C^*$-algebras, an ergodic action of a compact quantum group
  is not necessarily associated to a quantum subgroup  \cite{Podles_quantum_spheres}.
  
  This quantum homogeneous space was obtained in   \cite{PRergodic} as an application to rigid tensor $C^*$-categories of a Tannakian theorem for ergodic actions of compact quantum groups
  on unital $C^*$-algebras in the same paper,
  that characterizes properties of the spectral spaces associated to an ergodic action of a compact quantum group on a unital $C^*$-algebra. This Tannakian theorem in turn was the extension to the general case of a previous result by Bichon, De Rijdt,  and Vaes for ergodic actions of compact quantum groups for which the quantum multiplicity takes its maximum value
  given by the quantum dimension \cite{BDV}.
  
  In \cite{P_embedding}, \cite{PRergodic} we gave algebraic and categorical characterizations respectively of the cases where an ergodic action of a compact quantum group arises from a quantum subgroup.
  In \cite{P_rep} we characterized the compact quantum group ${\rm SU}_q(d)$ of Woronowicz in terms of its representation category.
  This work turned out to be closely related to an earlier work by Kazhdan and Wenzl for braided fusion categories of type $A$ fusion rules
  in the generic case \cite{KW}, with different methods, we built on previous work   by Doplicher and Roberts and Woronowicz.
  
In \cite{PRinduction}, we developed a quantum Mackey induction theory starting from a rigid tensor  $C^*$-category with simple tensor unit
and in particular in case of a generating object for the category, we associated a class of   compact quantum groups and  their ergodic $C^*$-actions   established in  \cite{PRergodic}. This   led us to  construct   natural finitely generated Hilbert bimodule representations of the acting  quantum group 
whose non-commutative $C^*$-algebra was our non-commutative homogeneous space. These bimodule representations
  form a
category  equivalent to the given category with all its structure. In the symmetric case, our construction reduced to the Doplicher-Roberts duality theorem with a different approach.
In this way we extended Doplicher-Roberts duality to  rigid tensor $C^*$-categories with simple tensor unit. 

By the generating property of the category, the compact quantum groups emerging from our construction were not unique but naturally associated to the category, identified with 
the universal compact quantum groups by Van Daele and Wang  \cite{W_VD}. We built on the classification of their irreducible representations by Banica  \cite{Banica1},  \cite{Banica2}.

 The next step remains the extension of the notion of compact quantum group to a suitable larger class that would admit natural examples with unitary braided symmetry. 
This step is necessary  to extend our   framework  to rigid, {\it unitarily} braided tensor $C^*$-categories, and to seek 
   applications in low-dimensional   algebraic quantum field theory.
The aim is the    construction of a quantum gauge group and, if necessary, a quantum homogeneous space on which the quantum group   acts naturally.

The first natural example of quantum gauge group within the braided Doplicher-Roberts program was given by Mack and Schomerus  
for the critical chiral Ising model of conformal field theory, described by 
the so called truncated fusion rules of the Verlinde category associated to ${\mathfrak sl}_2$ at the fourth root of unity \cite{MS_Conformal_field_algebras}.

In \cite{CP} we studied the structure of Mack and Schomerus example in detail and extended it   to the construction of an operator algebraic quantum group $A_W({\mathfrak sl}_N, q, \ell)$
over ${\mathfrak sl}_N$ that furthermore admits unitary representations of the braid group, and is naturally associated to the fusion category
${\mathcal C}({\mathfrak sl}_N, q, \ell)$ of the Drinfeld-Jimbo quantum group $U_q({\mathfrak sl}_N)$ at roots of unity. To this aim, we built on the work of Wenzl on the unitary structure of ${\mathcal C}({\mathfrak sl}_N, q, \ell)$ \cite{Wenzl}.
Moreover we gave a categorical characterization of the representation category of our examples $A_W({\mathfrak sl}_N, q, \ell)$
extending the results of \cite{P_rep}.

The inner structure of the examples of \cite{MS_Conformal_field_algebras} and \cite{CP} has been investigated in \cite{FKL}, and extended to all the Lie types $A_W({\mathfrak g}, q, \ell)$. A first particularly important property is
that of a {\it weak Hopf algebra} in a new sense as compared to the work by Mack and Schomerus and previous work in the literature.
Another relevant property is its {\it coboundary symmetry} associated to its $R$-matrix, since it leads to the desired extension of the notion of symmetric functor  in the category of its representations.

In \cite{FKL} we have extended Woronowicz theory beyond compact quantum groups to a new kind of weak Hopf $C^*$-algebras with sufficiently general unitary structures, among them a unitary structures induced by the coboundary symmetry. This specific case
is responsible of   unitary representations of the braid group. We extended the theory   to weak quasi-Hopf $C^*$-algebras with similar symmetry to accomodate examples associated to model of
unitary vertex operator algebras or conformal nets in a natural way, by unifying
with Drinfeld setting, using methods of unitary structures suggested by the work of Wenzl \cite{Wenzl}.
The earlier examples of \cite{MS_Conformal_field_algebras}, \cite{CP} and the new examples \cite{FKL} turn out to be   examples of such generalizations as quantum symmetry groups of the WZW model in the setting of vertex operator algebras. They play a primary role on all the main results of \cite{FKL}.

  In \cite{FKL}, we   investigated Doplicher-Roberts program   for some models of conformal field theory in the operator algebraic approach
  and the vertex operator algebra approach. We considered in particular   models with  pointed fusion categories and the WZW model in the context affine vertex operator algebras at positive integer levels. 
  Our work explores and establishes connections of the   braided Doplicher-Roberts program with a problem posed by Yi-Zhi Huang regarding a direct proof of the Finkelberg-Kazhdan-Lusztig theorem \cite{Huang2018}. 
  
  Due to what is known in the literature regarding the generators  of the quantum group  fusion categories ${\mathcal C}({\mathfrak g}, q, \ell)$
  for all the classical   Lie types and $G_2$ -- works  extending Schur-Weyl duality to   more Lie types with Brauer-Schur-Weyl duality and to the quantum case in the generic case and at roots of unity -- such  categories  may be regarded the braided analogues of the minimal categories in the    Doplicher--Roberts setting.  
  
  As a  consequence, this generating property also holds for their counterparts coming from conformal field theory
  such as ${\rm Rep}(V_{{\mathfrak g}_k})$ in the setting of the affine vertex operator algebras for the Finkelberg-Kazhdan-Lusztig equivalence, or the corresponding loop group conformal net due to Gui's isomorphism of categories between the vertex operator algebra categories and the corresponding conformal net categories \cite{Gui_isomorphism}.
  
 Specifically, in the quantum group case the braid group generating property in
 the Lie types $A$ and $C$ for the vector representation are well known and at the heart of the origin of quantum groups. If $q$ is not a root of unity, the  centralizer algebras 
$(V^{\otimes n}, V^{\otimes n})$ in the type $A$ case are described
by Hecke algebras, a quotient of the complex braid group algebra by well known generators and relations    Jimbo \cite{J}.
For the  Lie types $B$, $C$, $D$, when $V$ is the vector representation, this is the BMW algebra studied by
Murakami and Birman and Wenzl   \cite{Murakami}, \cite{Birman_Wenzl}. See also the exposition in \cite{GHJ}. 
Duality for the vector representation is of quantum groups was studied by Kirillov and Reshetikhin \cite{Kirillov_Reshetikhin}.
A summary of these results and references to the original papers for the vector representation for all the classical Lie types may be found in Theorem 10.2.5 of \cite{Chari_Pressley}.
The centralizer algebra of the fundamental representation for types $B$ and $D$ may be found in \cite{Wenzl_dualities} for the even case, and 
in \cite{Wenzl_dualities0} for the odd case. A special case has been studied by Rowell and Wenzl in 
\cite{Rowell_Wenzl}.
Type $G_2$ was shown by Lehrer and Zhang \cite{Lehrer_Zhang}, Morrison  \cite{Morrison}, and Martirosyan and Wenzl
 \cite{Martirosyan_Wenzl}.

  We take the generating property of the braid group representation between tensor powers of the fundamental representation as a synthesis of these examples within the braided Doplicher-Roberts program,
  in that it gives a conceptual explanation of the reason why these are examples of braided categories for which we can reconstruct
  quantum gauge group $A_W({\mathfrak g}, q, \ell)$ in a natural way directly from the category, with no need of a quantum analogue of Mackey induction
  for these models.
  
  As a matter of fact, we make use of this generating property to complete our  proof of the Finkelberg-Kazhdan-Lusztig theorem started in \cite{FKL} for the classical Lie types and $G_2$ with respect to Huang-Lepowsky braided tensor structure, see Sect. 11 in \cite{Pinzari_constructing}. Given this, we interpret Huang-Lepowsky braided tensor structure in natural connection with the braided Doplicher-Roberts program. 
  
  In   \cite{FKL} we   have
  constructed   a rigid ribbon braided tensor category structure on ${\rm Rep}(V_{{\mathfrak g}_k})$ equivalent to that of the quantum group fusion category via   the Drinfeld twist method. 
  For the Lie types $E$ and $F$ we do not know the generating property for  ${\mathcal C}({\mathfrak g}, q, \ell)$, therefore we can not conclude on complete equality of our   braided tensor structure with Huang-Lepowsky braided tensor structure on ${\rm Rep}(V_{{\mathfrak g}_k})$, but only partial equality on many braiding and associativity morphisms. 
    We shall discuss more the relation between the Finkelberg-Kazhdan-Lusztig theorem and the braided Doplicher-Roberts program from Sect. 3 on, especially on its role of unitarizing braided tensor structure of affine vertex operator algebras at positive integer level.

We approached the   problem   of the extension of the notion of symmetric functor by passing to the notion of coboundary,  first introduced
  by Drinfeld in the setting of his Drinfeld-Kohno theorem \cite{Drinfeld_quasi_hopf}. In Drinfeld setting, a coboundary is obtained from the braided symmetry associated to a Drinfeld-Jimbo quantum group, which admits a
    ribbon structure. 
    
    In \cite{FKL} we defined and found  natural coboundary structures in the fusion category of  quantum groups at suitable roots of unity 
    ${\mathcal C}({\mathfrak g}, q, \ell)$, and in the category of modules of an affine vertex operator algebra $V_{{\mathfrak g}_k}$ at a positive integer level that are induced by corresponding weak quasi-Hopf algebras, which are one the twist of the other,
    by a kind of twist that leaves a couboundary structure invariant. We relate these coboundary structures  
    to   the permutation symmetry via a functor. This is possible because the couboundary is not a representation of the braid group.
    This is our extension of the notion of symmetric tensor functor. 
    
    Slightly more  precisely, we have constructed an   suitable weak Hopf $C^*$-algebra $A_W({\mathfrak g}, q, \ell)$  naturally associated to  ${\mathcal C}({\mathfrak g}, q, \ell)$ and a
    related weak quasi-Hopf $C^*$-algebra structure on the Zhu algebra $A(V_{{\mathfrak g}_k})$ \cite{Zhu} associated to ${\rm Rep}(V_{{\mathfrak g}_k})$, and the two algebras are linked  by a direct explicit coboundary preserving relation.
    Our methods were suggested by the already mentioned work by Drinfeld on Drinfeld-Kohno theorem and   Wenzl on the unitary structure of  the fusion category of a compact quantum group at roots of unity.
    
    In the final section, we discuss some open problems and how the braided Doplicher-Roberts program might eventually address the challenges of 
 gauge theories and holographic correspondences.

     The material presented here originated as an introductory section of \cite{FKL}, which
     was omitted from the final version, and   subsequently expanded into the present   article to allow for a more detailed treatment.

\smallskip

\section{Quantum integrable systems, quantum groups, and  CFT}
  
\noindent{\bf  Baxter's work, the discovery of quantum groups and their appearance in  CFT.}  The study of completely integrable systems saw an increase of interest during the 70s, when powerful methods for the study of exactly solvable models were developed in the classical case with
the method of inverse scattering, and extended to the quantum case, see the review paper by Faddeev on the state of art till  1973 \cite{Faddeev1}, \cite{Faddeev2}. 
Complete integrability of models of conformal field theory   has been studied by many authors since the mid 70s.
The case of classical relativistically invariant  models was studied by Faddeev and his collaborators and led them to  investigate corresponding quantum versions to study quantum fields of completely solvable models.
The first attempts  in the quantum case were not satisfactory in that they led to   approximate deductions.   

In the early 70s Baxter's did breakthrough   work for certain models of statistical mechanics he was able to solve by introducing powerful methods \cite{Baxter}.
The quantum $R$-matrix and what later became to be  known as the Yang-Baxter equation was   established by Baxter for these models.  Baxter's work is responsible 
for the upcoming  discovery of {\it quantum groups}. See also \cite{Jones_in_and_around}, \cite{McCoy} for related survey papers.

In the late 70s Faddeev
and Sklynanin looked for a method that would unify the classical and the quantum cases,
and strikingly they found it with the use of the classical and quantum $r$- and $R$-matrices respectively.
The method of classical    $r$-matrix in the framework of the inverse scattering problem   was shown to admit a direct generalization to the quantum case. Quantum $R$-matrices were thus used   in the framework of the   quantum version of the inverse scattering problem and discovered symmetry structures given by quantum algebras see \cite{Sklyanin29}, \cite{FST},
\cite{Sklyanin}.

 The search
for solutions of the quantum Yang-Baxter equation has been of major interest   since Baxter work. 
%%The quantum $R$-matrix was studied to different settings.
In the mid 80s {\it quantum groups} were discovered as deformations of universal enveloping algebras by Drinfeld and Jimbo   \cite{D}, \cite{J} and named in this way by Drinfeld  \cite{Drinfeld_qg}. The term has soon been used far beyond the Drinfeld-Jimbo deformation setting.
Notably, the Yang-Baxter  equation is satisfied by the universal $R$-matrix of the Drinfeld-Jimbo quantum groups.
%%, to include the setting of earlier mentioned quantum algebras.

Earlier models in conformal field theory would in retrospect become the first models with symmetries given by quantum groups.
Several other quantum mechanical models were found to admit signs of this new kind of symmetry, such as \cite{AG}, \cite{Moore_Seiberg2},
\cite{Moore_Reshetikhin} to name just a few.
See also the introductions and references in \cite{MS}, \cite{Schomerus_field_algebra}, \cite{Chari_Pressley}, \cite{FKL}.
These early models based on the quantum $R$-matrix, may be regarded as the place where   the problem of statistics in a sense of extending the usual Bose-Fermi alternative, made his first appearance. This notion was greatly developed in the setting of algebraic CFT by several authors  in low dimensions, and we briefly next recall this.

\smallskip

  \noindent{\bf Statistics in algebraic quantum field theory, Doplicher-Roberts theory and the problem of extension to conformal field theory.} 
  In the algebraic approach to quantum field theory by Doplicher-Haag-Roberts, superselection sectors of the Hilbert space of physical states are described by representations of the observable operator algebras ${\mathcal A}({\mathcal O})$ satisfying localization principles motivated by relativistic quantum theory. The vacuum sector is recovered by a GNS construction
  from the algebra and the relevant representations   turn out to correspond to the localized endomorphisms of the   algebra. 
  One of the successes of algebraic quantum field theory is that existence of antiparticles  emerges naturally from the structure of superselection sectors, see \cite{DHR}, see also \cite{Haag}. One derives existence of creation and annihilation operators ${\overline R}$ and $R$ that describe for each localized endomorphism $\rho$ a conjugate, or dual, endomorphism $\overline{\rho}$.
 The localized endomorphisms form the structure of a symmetric rigid strict tensor $C^*$-category   in $4D$ space-time dimension. The tensor product is given by composition, and the permutation symmetry of the category arises from the permutation group statistics \cite{DHR}.
  
 Doplicher and Roberts  reconstructed  field algebras 
 ${\mathcal F}({\mathcal O})$  generated by  localized field operators, together with a compact gauge group $G$, a symmetry group of the theory,
which acts as a group of automorphisms of ${\mathcal F}$, such that the $G$-invariant elements
    are precisely the local observables algebras   ${\mathcal A}({\mathcal O})$  \cite{DR_CMP}.

It has been known for a long time that when $D\leq3$, the corresponding category is braided rather than symmetric.  The braided symmetry arises from the braid group statistics, already for endomorphisms localised in bounded regions for $2D$ and in unbounded cones    for $3D$  \cite{FRS}, \cite{F}, \cite{Frolich}.

Following the Doplicher-Roberts  reconstruction of the gauge group and field algebra from the net of local observable algebras and its localized endomorphisms endowed with permutation statistics,  there arose the question of whether their theory could be extended to algebraic conformal field theory,
with the gauge group replaced by a quantum group, and braid group statistics corresponding to a quantum $R$-matrix.
This problem gave rise to a series of works, see the introduction of \cite{FKL} and references therein.
Among these works, several authors addressed this problem starting from specific models.

 \smallskip

\noindent{\bf  Quasi-Hopf algebras.}
 In   1989, notably    Drinfeld introduced {\it quasi-Hopf algebras} as a   generalization of Hopf algebras,    relaxing the requirement of coassociativity through the introduction of a non-trivial associator. This put emphasis on non-strict braided tensor categories arising from models of conformal field theory.

Drinfeld constructed a braided quasi-bialgebra structure over $U(\mathfrak{g})$ with   braiding and associator induced by the   Knizhnik-Zamolodchikov (KZ)  equation.
 The resulting representation  category of $U(\mathfrak{g})$,    the Drinfeld category,   inherits the structure of a  braided tensor category.

 Furthermore, Drinfeld introduced the notion of {\it twist}, a mechanism that preserves a quasi-Hopf algebra structure. 
 In the {\it Drinfeld-Kohno theorem}, Drinfeld constructed a twist that realizes a braided tensor equivalence between
  the representation category of the ribbon quantum group $U_{h}(\mathfrak{g})$ and the Drinfeld category, 
  establishing the latter as a rigid ribbon category \cite{Drinfeld_quasi_hopf}.

   Drinfeld commented on these constructions as approximate physical symmetries for the
  { WZW model } in conformal field theory (CFT) \cite{Drinfeld_cocommutative}. 
  Complex quasi-Hopf algebras over discrete group algebras later appeared in the work of Dijkgraaf, 
  Pasquier, and Roche on orbifold models  \cite{DPR}.
  \smallskip

 \noindent{\bf Weak quasi-Hopf algebras in CFT and unitarity.} In 1990, Mack and Schomerus first constructed 
a field algebra in the framework of Doplicher-Haag-Roberts theory,
associated to a model of conformal field theory, the critical chiral Ising model. In their example  the symmetry group $G$ is replaced by a semisimple quotient of   Drinfeld-Jimbo quantum group  $U_q({\mathfrak sl}_2)$, with deformation parameter $q$ given by a primitive fourth root of unity, and the Bose-Fermi local commutation relations of fields are replaced by  {\it local braid relations,} first proposed by Fr\"olich \cite{Frolich_local_braid}, see also  \cite{MS_Conformal_field_algebras} and references therein. This was the first example where such a generalization of Bose-Fermi commutation relations was obtained in low dimensional quantum field theory. Thus the example was interesting and suggested the question of what kind of quantum group they obtained.
In their case, they needed to exclude many representations of the quantum group
considered as unphysical, and retained only finitely many \cite{MS_Conformal_field_algebras}.

In subsequent papers, Mack and Schomerus studied their quantum symmetry group further,
and, motivated by Drinfeld work on quasi-Hopf algebras,  proposed a theory of quantum symmetries via  {\it semisimple weak quasi-Hopf algebras}   within the framework of   conformal quantum field theory
\cite{MS_Conformal_field_algebras}. In particular, their previous first example   being a natural semisimple quotient of $U_q({\mathfrak sl}_2)$ as an algebra, was defined as endowed with a natural new coalgebra structure compatible with the semisimple fusion rules of the category of localized endomorphisms in  \cite{MS1}, \cite{MS}. In this way they used weak coproduct to describe   semisimple {\it truncated tensor products} 
from the non-semsimple quantum group $U_q({\mathfrak sl}_2)$, and regarded their weak quasi-Hopf algebra as the {\it true symmetry of the minimal conformal models}.
In  \cite{MS} Mack and Schomerus recognized the value of their constructions as {\it the general algebraic structure of a symmetry in low dimension}. 

Their constructions benefitted from the special fact that tensor products of irreducible representations of ${\mathfrak sl}_2$ is multiplicity-free. It is therefore desirable to go beyond the case of ${\mathfrak sl}_2$.

 Moreover, while they defined weak quasi-Hopf algebras   with non-unital coproducts   -- allowing for an analogue of twist deformation -- their framework omitted a corresponding {\it semisimple weak Hopf algebra structure}. Identifying such a structure is essential for recovering Drinfeld's 3-coboundary associator achieved by the Drinfeld-Kohno theorem. 
Additionally,    the unitary structure of tensor products of their weak quasi-Hopf algebras was overlooked in their work, and later studied in a close setting
by Wenzl \cite{Wenzl} and Xu \cite{Xu_star}, whose results set up an earlier conjecture by Kirillov on positivity of an hermitian form
associated to the fusion category of quantum groups at roots of unity.
Wenzl's approach opened up  to new light on the mentioned problems, and are at the base of our previous work \cite{CP}, \cite{FKL}.

\smallskip

   \noindent{\bf Quantum symmetries in CFT. }
   Thus, there remains the   problem of unifying the two approaches: the  classical gauge group case by Doplicher and Roberts in the algebraic approach to  $4D$  QFT
   and the quasi-Hopf or weak quasi-Hopf algebra case in  CFT, in the approaches by Drinfeld and Mack and Schomerus, respectively, into one method of investigation of quantum fields. 

        The need of a natural notion of   weak Hopf algebra among the semisimple weak quasi-Hopf algebras is twofold in the setting of conformal field theory.
   On one side, it naturally emerges from the desire of extending and reconciling with
   the construction by Doplicher-Roberts of  a  compact gauge group and the field algebra in  $4D$.    Moreover, in their case, the tensor category is strict, and the compact group is unique up to isomorphism.
   On the other, while the first example by Mack and Schomerus of weak quasi-Hopf algebras in the case of the Ising model is very naturally associated to the model,
   and has the   virtue of allowing the field algebra  construction as shown by  Schomerus in \cite{Schomerus_field_algebra}.
   This construction is not entirely satisfying due to  the use of   non-trivial associators brought complications to this construction, and one would like to look for a better understanding of the involved associators.    In this  case, one also needs to look for canonical association of weak quasi-Hopf algebras.

   Moreover, the mentioned field algebra construction   is extended to the general case and holds for any rational theory, but in these more general
    cases, the associated weak quasi-Hopf algebras are constructed using only a weak dimension function associated to the category.
     This
   leads to high non uniqueness and appearance of non-canonical weak quasi-Hopf algebras as compared to the case on high dimension
   and also to the examples of the minimal models treated in \cite{MS_Conformal_field_algebras}, \cite{MS1}.
   
   In the  Doplicher-Roberts case 
     the advantage of uniqueness allows to overlook  at concrete   properties of relativistic quantum field theory   that might prove useful to overcome lack of uniqeness in low dimensions, and this poses the problem of finding an efficient replacement in low dimensions where uniqueness of the functor from the category of localized endomorphisms to the category of vector spaces does not hold, and we would like to look for canonical natural constructions.
   
   Our proposal is that of using the {\it minimum energy functor} of affine vertex operator algebras, which exists in all vertex operator algebras,  to construct such a natural
   quantum gauge group, that turns out to be a Drinfeld twist of a weak Hopf algebra in a new sense defined in \cite{FKL}, with a braid-related unitary structure, in a way that
   follows the suggestions of the Drinfeld-Kohno theorem but in the setting of the axiomatic approach  to conformal field theory provided by vertex operator algebras.
   
   As mentioned in Sect. 1, for a long time, the compact quantum groups by Woronowicz  
 were considered quite far from applications to low dimensional algebraic quantum field theory
 because of lack of examples with unitary braided symmetry.

Our previous work  \cite{FKL} extends Woronowicz's theory beyond compact quantum groups to accommodate this class of representations. By unifying his framework with Drinfeld's setting and employing the unitary structure methods suggested by Wenzl [53], we establish the main results of \cite{FKL}   based  on canonical constructions of such generalizations as quantum symmetry groups of the WZW model in the setting of vertex operator algebras.
   \smallskip

\section{The Finkelberg-Kazhdan-Lusztig equivalence} 

These problems naturally led us to the problem posed by Huang, mentioned at the beginning.
Fundamental results by Kazhdan, Lusztig, and Finkelberg establish a rigid braided tensor category structure on the category $\tilde{\mathcal O}_\ell$ of modules of an affine Lie algebra at a positive integer level and its braided tensor equivalence with a   semisimple subquotient fusion category
${\mathcal C}({\mathfrak g}, q, \ell)$ of  a quantum group at roots of unity $U_q({\mathfrak g})$ \cite{KLseries}, \cite{Finkelberg}, \cite{Finkelberg_erratum}.
For an overview   of these results, see e.g.
Sect. 1 of \cite{FKL}, see also \cite{Huang2018} or the introduction of 
\cite{Pinzari_constructing}.

The original proof   is considered   indirect, as  the equivalence   between $\tilde{\mathcal O}_\ell$ and ${\mathcal C}({\mathfrak g}, q, \ell)$
is given by the composition of two separate equivalences.  

While the Kazhdan-Lusztig equivalence naturally relates two non-semisimple representation categories defined by   negative
rational levels of affine Lie algebras and root of unity deformations of universal enveloping algebras, from which semisimple structures are obtained as subquotients in the exposition by Finkelberg, the necessary  passage to the negative levels of the composed equivalence is often viewed as unnnatural.

   Furthermore, Finkelberg's proof of rigidity for $\tilde{\mathcal{O}}_{l}$
    requires the rigorous establishment of the Verlinde formula and the computation of fusion rules,
     involving extensive work by Faltings and Teleman. 
     Huang had earlier first noted that the Verlinde formula is essential for establishing rigidity and modularity within the representation category of vertex operator algebras (VOAs). 
Thus the construction of   rigid braided tensor category structures  on   categories of certain modules on the
affine Lie algebra or vertex operator algebra at positive integer levels, becomes substantially more difficult to follow as compared to the conceptual clarity offered by the Drinfeld-Kohno theorem.

See   \cite{Huang2018}, the introduction of \cite{Pinzari_constructing}, Sect. 1 in  \cite{FKL},
 for more details.
For an overview   of these results, see e.g.
Sect. 1 of \cite{FKL}, see also \cite{Huang2018} or the introduction of 
\cite{Pinzari_constructing}.
For more background material, the reader may consult   results in algebraic quantum field theory, vertex operator algebras, quantum groups, and tensor categories as detailed in Section 1 of \cite{FKL} and the references therein.

     These complexities motivate a search for a more direct, self-contained proof operating as proposed by Huang in  \cite{Huang2018},
     within a semisimple framework, that hopefully simplifies the proof of rigidity.

 \smallskip

\noindent{\bf Our strategy for a direct proof}\label{2}
Our work stems from the observation that the original proof of the Kazhdan-Lusztig-Finkelberg theorem does not utilize 
semisimple weak quasi-Hopf algebras, because the involved categories lack an initial
 fibre functor, unlike those in the Drinfeld-Kohno theorem.
 Therefore this links   to an abstract duality problem extending the work of Doplicher and Roberts on their
abstract duality theory for compact groups.

In the paper \cite{FKL} we extend Drinfeld's approach to the Drinfeld category by introducing weak Hopf algebras as quantum symmetries for the braided fusion categories  associated to the affine vertex operator algebras at positive integer levels. In our model, weak Hopf algebras and weak quasi-Hopf algebras occupy positions analogous to Hopf and quasi-Hopf algebras in original Drinfeld theory.

 To prove the Finkelberg-Kazhdan-Lusztig Theorem in the context of vertex operator algebras, we identify a semisimple quotient algebra, 
 $A_{W}(\mathfrak{g},q,\ell )$, of the quantum group $U_{q}(\mathfrak{g})$ at a root of unity. This quotient plays a role similar to the formal deformation $U_{h}(\mathfrak{g})$ in the Drinfeld-Kohno theorem, while the Zhu algebra $A(V_{\mathfrak{g}_{k}})$ associated to
 the vertex operator algebras serves as the analogue to the quasi-Hopf algebra $U(\mathfrak{g})[[h]]$. 
 
 The most challenging part is how to define a coproduct, and how to put on  $A_{W}(\mathfrak{g},q,\ell )$ all the structure
 that resembles the structure that Drinfeld put on $U_{h}(\mathfrak{g})$ to prove the Drinfeld-Kohno theorem.
 A main difficulty is that $A_{W}(\mathfrak{g},q,\ell )$ does not result from a coideal of $U_{q}(\mathfrak{g})$ because
 the coproduct of $U_{q}(\mathfrak{g})$ is not compatible with the tensor product of the fusion category ${\mathcal C}({\mathfrak g}, q, \ell)$.
 Thus these constructions present technical difficulties, and the fact that all the properties can be established gives an astonishing
 meaning to Drinfeld work in the Drinfeld category.
 
    Our approach utilizes Wenzl's unitary structure analysis of \cite{Wenzl} in an essential way, which is not present in Drinfeld work,
   to endow  $A_{W}(\mathfrak{g},q,\ell )$ with a $C^*$-structure
   compatible with the unitary  ribbon-braided tensor structure of the quantum group fusion category, that we call {\it unitary coboundary}, extending both Drinfeld's setting    and Doplicher-Roberts setting for the symmetric functor associated to the representation category of
 a classical compact group.

Then the unitary structure allows the construction of a Drinfeld twist
analogous to the original Drinfeld twist in the algebraic setting, and is useful  to transfer all the   structure, establishing     $A(V_{\mathfrak{g}_{k}})$ as  unitary coboundary weak  quasi-Hopf $C^{*}$-algebra. In particular, the twist operation induces unitary structure, tensor product, braiding and ribbon structure, and rigidity on the Zhu algebra. In this way the representation category of the
vertex operator algebra becomes enriched with all this structure thanks to the well-known Zhu linear equivalence and its inverse.

 Thus by equipping $A_{W}(\mathfrak{g},q,\ell )$  with  the compatible unitary coboundary weak Hopf algebra structure from quantum group structure, we 
 derive  a unitary coboundary weak quasi-bialgebra structure on  $A(V_{\mathfrak{g}_{k}})$,
 by  establishing a Drinfeld twist that relates their respective fusion categories. This twist allows to prove    equivalence of the braided tensor structures, the original one and the twisted one,
 and characterizes the weak Hopf algebra $A_{W}(\mathfrak{g},q,\ell )$ as a compact group-like object.
 
 Then we prove that,   for the classical Lie types and $G_2$, the   ribbon-braided tensor structure obtained on the representation catgeory of the vertex operator algebra identifies with the structure independently introduced by Huang and Lepowsky 
 in a general framework. 
 
 To this aim  we rely on some keypoints. The first one is the identification of the fundamental notion of primary field that hints at the decomposition of tensor products of two representations of the underlying compact group following the Verlinde fusion rules, by focusing on the minimum energy functor in the setting of affine
 Lie algebras (which corresponds to Zhu's functor in the framework of VOAs). This analysis gives lights to the   unitary structure,   the braided symmetry, the fusion rules on special tensor product of objects that become coinciding with those known for affine VOAs.
 
 These special tensor products are made of   a fundamental representation of the compact group as one of the two factors, and an arbitrary irreducible as the other factor. These space of intertwiners are tied to corresponding spaces of intertwiners of fusion rules of the
 quantum group, by the unitary structure. Thus the Verlinde fusion rules are tied to the fusion of the quantum group in a canonical way by the unitary structure.
 
 This construction relies on the work by Wenzl for all Lie types on the quantum group side, and
 Wassermann, and Toledano-Laredo on the affine Lie algebra side.  Moreover, our twist that makes the unitary structure of the quantum group weak Hopf algebra into a trivial one, also shows identification of  
the ribbon, braiding, and associators at the level of weak quasi-bialgebra structure on the Zhu algebra, inspired by methods of
Drinfeld. 

Finally, we
 develop
  a uniqueness theorem
 for the associator in case of two braided symmetries (the one introduced independently in the framework of loop groups and Huang-Lepowsky work and the other
 introduced by our twist operation from quantum groups) that coincide on special pairs of morphisms.
 For the Lie types $E$ and $F$ we prove coincidence of the braiding morphisms (and associativity morphisms) on many pairs (or triples) of objects using a generating representation as one (or two) of the factors, see Theorem 2.4 in \cite{FKL}.

  \medskip

\section{Manifestly unitary structures in models of CFT} The notion of unitarity plays a central role in \cite{FKL}. 
Inspired by Vaughan Jones's insights into CFT, one expects the structures under consideration to display manifest unitarity. However, Wenzl's work on quantum groups at roots of unity complicates this view by showing that the unitary structure of the associated fusion categories is twisted on tensor products \cite{Wenzl}. While this contrasts with the standard operator algebraic approach to quantum groups -- which typically employs trivial tensor products of unitary structures -- and also calls its manifest unitarity nature into question,  this twist is precisely what enables unitary representations of the braid group. Since braid group unitarity is a fundamental feature of both conformal nets and the models studied in that setting, a   question arises: how can these two seemingly disparate aspects of unitarity of tensor products be reconciled?
\smallskip

In the setting of    loop group conformal nets  manifest unitarity
may be observed for the WZW model in the type $A$ case,  they correspond  to certain  primary fields   that appear in works   by
     Wassermann and    Toledano Laredo see Sect. IV in \cite{Wassermann} and \cite{Toledano_laredo}. 
     At the level of the top space module of the simple Lie algebra, the fusion tensor product appears manifestly unitary. 
This was the    starting point of our work in constructing the connection.

 Gui extended the work of Wassermann and Toledano-Laredo  to the setting of VOAs. He showed unitarizability and tensor equivalence with the representation category of the corresponding conformal net \cite{GuiI}, \cite{GuiII}, \cite{GuiIII}, \cite{Gui_isomorphism}.\medskip

\section{Algebraic, Hermitian and unitary theory of weak quasi-Hopf algebras.}

  H\"aring-Oldenburg extended to the weak case the work of Majid \cite{Majid2}, \cite{Majid4} for quasi-Hopf algebras \cite{HO}.
He defined a weak quasi-tensor fibre functor from a semisimple rigid tensor category and formulated a  Tannaka-Krein duality theorem showing  that such pairs   are in duality with semisimple weak quasi-Hopf algebras. This duality   relates a braiding of the category to a  quasitriangular structure of the algebra.

It follows from \cite{HO} that  a necessary and sufficient condition for  ${\mathcal F}$ 
to be upgraded to a weak quasi-tensor functor is that   $\rho\to {\rm        dim}({\mathcal F}(\rho))$ be a weak dimension function, meaning that ${\rm        dim}({\mathcal F}(\rho\otimes\sigma))\leq {\rm        dim}({\mathcal F}(\rho)){\rm        dim}({\mathcal F}(\sigma))$ for all irreducifble objects $\rho$, $\sigma$. 
A weak quasi-tensor structure with the same dimension function  is  not unique,   passing to another affects   the weak quasi-Hopf algebra by a twist deformation. 
More than this, there may be infinitely many weak dimension functions. These matters of high non-uniqueness 
are addressed in \cite{HO}, \cite{Schomerus_field_algebra}, and contrast with the case of Doplicher-Roberts construction of the compact group.

For this reason, the approach of \cite{FKL}  is to start with  a semisimple tensor category ${\mathcal C}$ possibly with more structure, endowed with a given linear      functor ${\mathcal F}:{\mathcal C}\to{\rm        Vec}$, that we understand  as naturally associated to ${\mathcal C}$, and soon we shall specify ${\mathcal F}$.

We discuss three main istances, {\it Wenzl functor} for the fusion category ${\mathcal C}({\mathfrak g}, q, \ell)$ of quantum groups at roots of unity,   {\it Zhu's functor} for the fusion category of a vertex operator algebra
 or finally {\it the minimum energy functor} in categories arising from models described in the framework of algebraic conformal field theory. Moreover, we look for naturally associated weak (quasi)-Hopf algebras with coproduct compatible with the tensor structure of the category.
 \smallskip

 \noindent{\bf Algebraic structure.} We develop a   theory of weak quasi-Hopf algebras over the field
of complex numbers and the Tannakian formalism between them and tensor categories.
We discuss extra structure such as quasitriangular and ribbon structures corresponding categorically
to a braiding and a ribbon structure. 

Then we  make a proposal of the weak analogue of Hopf algebras, that we call 
 weak  Hopf algebras,    among weak quasi-Hopf algebras  in a cohomological interpretation.    We also  introduce a notion of weak tensor functor between tensor categories. A weak  Hopf algebra is   characterised, via Tannaka-Krein duality,   by a semisimple rigid tensor category endowed with a weak tensor   functor to ${\rm        Vec}$.

 We develop a theory for weak  Hopf algebras which includes  the notion of   $2$-cocycle deformation,   quasi-triangular and ribbon structure. We introduce  twisted Hermitian or $C^*$-structures, and study the relationship with the ribbon structure, and with  unitary braided symmetry and coboundary symmetry for the representation category.  In particular, we introduce the notion of unitary coboundary weak  Hopf $C^*$-algebra.  
We show that the examples $A_W({\mathfrak sl}_N, q, \ell)$ associated to the fusion category of $U_q({\mathfrak sl}_N)$ at roots of unity as developed in  a previous paper \cite{CP} and the further examples  $A_W({\mathfrak g}, q, \ell)$  corresponding to  $U_q({\mathfrak g})$ constructed in \cite{FKL}  are of this kind.\smallskip

\noindent{\bf Hermitian and unitary structure.} 
 In the framework of Tannakian duality, the unitary structure of a unitary tensor category manifests within the associated weak quasi-bialgebra as a positive $\Omega $-involution.

   The operator $\Omega $ is a positive element in the tensor square of the algebra that measures the failure of the $*$-involution to 
 commute with the coproduct.  
  Specifically, $\Omega $ becomes trivial if the involution and coproduct commute -- a property we expect only in cases of manifest unitarity.

For a general $\Omega$-involutive weak quasi-bialgebra, the category of $^*$-representations on Hermitian spaces 
forms a   tensor $^*$-category. In this setting,  the   hermitian form of the tensor product
is defined by  the action of $\Omega$. When $\Omega$ is positive, the full subcategory of 
Hilbert space representations is a unitary tensor  category. 

Furthermore, we establish the rigidity property
across all three settings provided an antipode exists. 
For unitary discrete  weak  Hopf $C^*$-algebras,
we provide an explicit construction of conjugates, extending a result by Woronowicz \cite{Wor} for standard  $C^*$-involutions.

Finally, by utilizing these unitary weak quasi-Hopf structures, we prove the uniqueness of unitary structures 
across a broad class of tensor categories. This result solves a problem posed by Galindo \cite{Gal}. 
Closely related results recently obtained by Reutter through alternative methods  \cite{Reutter}.
 
 The instances arising from quantum groups at roots of unity are described by specific $\Omega$-involutions related to
 the braiding structure, and we next recall this.
 A striking discovery by Wenzl reveals that, for the non-semisimple quantum groups at roots of unity, the non-trivial operator $\Omega $ coincides precisely with the Drinfeld coboundary matrix $\overline{R}$ -- the same matrix Drinfeld employed to prove the Drinfeld-Kohno theorem for the Drinfeld category. This correspondence suggests a profound connection between unitary structures and tensor equivalences, a relationship we explore  some more in this work and in detail in \cite{FKL}.
 \smallskip

\noindent{\bf  The striking interplay between the braiding and unitary structures.}
The theory of Hermitian (unitary) coboundary wqh algebras given in \cite{FKL} is a special case of $\Omega$-involution. It has a twofold
motivation in \cite{FKL}.    
On one side it goes back to Drinfeld ideas of coboundary matrix in the setting of quasi-Hopf algebras (a $2$-coboundary deformation of the $R$-matrix)
which plays a role in the proof
of the original Drinfeld-Kohno theorem. On the other side, it follows closely the study of unitarity of the fusion categories
associated to quantum groups at roots of unity by Wenzl. For example, in our terminology,
Wenzl showed, among other things, that Lusztig integral form of $U_q({\mathfrak g})$ at roots of unity is an Hermitian coboundary Hopf algebra over a suitable polynomial ring in  a formal variable
and used this structure to show a conjecture of Kirillov about positivity of a  $^*$-structure on the associated fusion ${\mathcal C}({\mathfrak g}, q, \ell)$ category for certain "minimal" roots of unity $q$. Similar results were obtained by Xu with different methods \cite{Xu_star}.

The unitary coboundary structure of the integral form  motivates us to look for 
non-formal semisimple unitary coboundary weak Hopf structures associated to fusion categories ${\mathcal C}({\mathfrak g}, q, \ell)$ of quantum groups at roots of unity, and we shall come back to this.\medskip

 \section{The interplay between braided tensor equivalences and unitary structures}

\noindent{\bf   Deriving unitary tensor structures from tensor equivalences  
} In the conformal net approach to CFT, or in the study of quantum groups at specific roots of unity,
the associated braided tensor categories are unitary.
Therefore, establishing a tensor equivalence with
   with another tensor category,  requires first to unitarize the latter, as a necessary condition.

 If our target category is concrete, such as the module category of a VOA, a natural first step
 is constructing   a unitary structure on its simple modules. We may do this guided by an explicit  linear equivalence
  between
   the examples under study: the explict form of 
  the linear equivalence will give insight on how to   unitarize the simple modules of the target category. 
  Once this linear unitarization is achieved, the next task
  is that of making the linear unitary structure of the target category into a unitary tensor structure.

Abstracting the situation,  if ${\mathcal C}$ denotes the tensor category and
we know  that  ${\mathcal C}$ is tensor equivalent to a unitary tensor category ${\mathcal D}$ and we   know how it acts as a linear 
equivalence, then we derive a linear $C^*$-category ${\mathcal C}^+$ linearly equivalent to ${\mathcal C}$.  
At this point   we may make ${\mathcal C}^+$ into a unitary tensor category via this equivalence. 
We discuss this abstract construction in Sect. 15 of \cite{FKL} with the help of weak quasi-Hopf algebras of Mack and Schomerus.

The following is a main example. Let    ${\mathcal C}$ be the category $\tilde{\mathcal O}_k$
associated to an affine Lie algebra at a positive integer level $k$ as considered by Finkelberg \cite{Finkelberg}, \cite{Finkelberg_erratum},
and for ${\mathcal D}$ we take the fusion category  ${\mathcal C}({\mathfrak g}. q, \ell)$ of the associated quantum group,
 which is unitary by the results of Wenzl and Xu  \cite{Wenzl}, \cite{Xu_star}, and on Finkelberg-Kazhdan-Lusztig theorem for the equivalence between ${\mathcal C}$ and ${\mathcal D}$.
 
 This construction shows that there will be two unitary weak quasi-bialgebras associated to the two categories,
 which are one the twist of the other, since the two categories are tensor equivalent.
 
One   then needs to understand the starting tensor equivalence in a direct way, 
to provide explicit constructions of unitary tensor structures on ${\mathcal C}^+$ with our method.\medskip

 \noindent{\bf Reconstructing tensor equivalences from unitary structures} 
%%Having established how tensor equivalences induce unitary structures, we now examine the converse relationship.
If we reverse the argument and rely on the idea that unitary structures in CFT are manifestly unitary on tensor products,
 then we may try to construct a tensor equivalence between representation categories of affine VOAs and quantum groups via a construction that makes the unitary structure of quantum groups trivial.

 The weak quasi-Hopf algebras   will   be useful  
 to {\it construct} tensor equivalences between tensor categories, and prove the Finkelberg-Kazhdan-Lusztig theorem in the setting of modules of vertex operator algebras with Huang-Lepowsky tensor structure, by looking for a trivial   unitary structure transported from the non-trivial unitary structure of tensor products of representations of quantum groups discovered in  \cite{Wenzl}.   
 
 It is quite striking that  unitary structures for tensor product representations are deeply linked to
 the equivalence between categories associated to conformal field theory with the corresponding Knizhnik-Zamolodchikov differential equation.

   Frenkel and Zhu \cite{Frenkel_Zhu},  
   associated a simple vertex operator algebra $L_{k, 0}$  
 to an affine Lie algebra $\hat{\mathfrak g}$. This   is a simple rational $C_2$-cofinite vertex operator algebra for $k\in{\mathbb N}$
 by  \cite{DongLiMason}. The simple  $L_{k, 0}$-modules $L_{k, \lambda}$ are classified by a finite set, the Weyl alcove
 $\lambda\in\Lambda^+_k$. By a fundamental result of Kac, Theorem 11.7 in \cite{Kac2}, $L_{k, \lambda}$ is unitarizable.
 By a result by
 Dong and Lin \cite{DongLin}  $L_{k, 0}$ is a unitary vertex operator algebra and $L_{k, \lambda}$ form a complete set of unitary
 $L_{k, 0}$-modules for $k\in\Lambda_k^+$.

 Sect. 21 of \cite{FKL} is introductory to vertex operator algebras, while Sect. 22 deals with
  various  examples  of VOAs, including the affine VOAs.  We give general results characterizing the enrichment of  a unitary braided tensor category structure to the Huang-Lepowsky braided tensor structure,
  in terms of a corresponding unitary weak bialgebra structure on thes Zhu algebra, and then we construct such a structure in the affine cases
  in a canonical way via a Drinfeld twist method. In this way, the unitary structure of the Zhu algebra is   a trivial one (coproduct and involution commute) by the twist from our weak Hopf algebra $A_W({\mathfrak g}, q, \ell)$, and it realizes an equivalence of the module category of the vertex operator algebra
  with the quantum group fusion category, by construction. We shall come back to the twist   in Sect. \ref{7}, and to $A_W({\mathfrak g}, q, \ell)$ in Sect. \ref{8}.

 \medskip

\section{A twist making  the Zhu algebra $A(V_{\mathfrak{g}_{k}})$   into a manifestly unitary coboundary weak quasi-Hopf algebra; main Theorems   in \cite{FKL}} 
\label{7} An   aspect making the     $C^*$-case of   interest is that there are cases   where
$\Omega$ admits a square root twist, that is a twist $T$ such that $\Omega=T^*T$, $\Omega^{-1}=T^{-1}(T^{-1})^*$.  
In this way  the $\Omega$-involution of a unitary weak quasi-Hopf    algebra can be twisted into one   in the usual   sense
that is the $^*$-involution commutes with the coproduct, and this seems to connect with Jones idea.
While this square root construction always exists for $\Omega$-involutive quasi-Hopf $C^*$-algebras with $\Omega$ positive, it is not clear whether the same holds in the weak case with useful properties.  One of the main result of \cite{FKL} is the construction of a canonical square root in the weak $C^*$-case for the unitary structure $\Omega$ of the weak Hopf algebra $A_W({\mathfrak g}, q, \ell)$. It has the nice property of making   the
braiding in a very simple "exponential" form, that is a close connection with the original  Drinfeld-Kohno theorem.
\medskip

  In the application, the situation is more complicated than how described in the previous paragraph.
We need   a   weaker notion of square root of $\Omega$, and
 the precise definition is in our abstract version of the original Drinfeld-Kohno theorem, see Theorem 29.4 in \cite{FKL}.
The difficulty of finding a square root $T$, $T^{-1}$ as stated above   is related, in our main application of fusion category of quantum groups, to certain eigenvalues arising from the braiding (the non-formal Drinfeld coboundary matrix). 
In part b) of the same theorem we introduce a pointwise notion of square root, that is a square root for a given fusion tensor product  $\rho {\otimes} \sigma$ 
on which $\Omega$ acts.

 In the application of quantum groups at roots of unity it turns out that $\Omega$ admits
a square root in this sense when $\rho=V_\lambda$ is any irreducible of the fusion category and $\sigma=V$
is a fundamental representation. This construction builds on the results of \cite{Wenzl}. Correspondingly, in the CFT this corresponds to constructing a {\it manifestly unitary
structure for the primary fields} in the sense used by Wassermann, Toledano-Laredo, Gui in  \cite{Wassermann}, \cite{Toledano_laredo}, \cite{GuiII}. 

To do this, we need to extend Wenzl Hermitian form $(u, v)$ from the building blocks $V_\lambda\otimes V$ to general tensor products of modules in the Weyl alcove, by taking   a canonical square root of the action of the quantum Casimir operator, that in the terminology of quantum groups is the ribbon element.
This leads to a positive solution of the existence of Wenzl Hermitian form on the tensor product of any pair of irreducible representations of $U_q({\mathfrak g})$ in the closed Weyl alcove, which, however, can not   always be positive by non-semisimplicity.
For the building blocks   this  Hermitian form is positive for all Lie types, and we use the work of \cite{Wenzl}.
This gives a positive $\Omega$-involution on the weak  Hopf algebra. Then we show that a positive square root   holds locally on the building blocks. This is summarised by our main   Theorem 2.2 in \cite{FKL}.
Theorem 2.4 of the same paper addresses the comparison of our ribbon braided tensor structure with that introduced by Huang and Lepowsky,
and we find that the two structures are the same for the clasical Lie types and $G_2$.  The final part of the proof is contained in \cite{Pinzari_constructing} and relies on  
the generating property of the braid group representation between   tensor powers of the generating representation in the fusion category, which is known for these Lie types. Moreover to this aim we develop a uniqueness theorem for pairs constituted by a braided symmetry and an associator in a semisimple linear category with a tensor product. In this way we overcome the use of the KZ equation, which is not natural
in our setting. Notice that Drinfeld's original proof of Drinfeld-Kohno theorem also relies on a uniqueness theorem for associators
for Drinfeld category, but with a very different approach \cite{Drinfeld_quasi_hopf}.
For the remaining Lie types $E$ and $F$   Theorem 2.4 in \cite{FKL} gives coincidence of all the structure on many objects and all morphisms.

\medskip

\section{Unitary coboundary weak Hopf structure of $A_W({\mathfrak g}, q, \ell)$}\label{8}

We found   in \cite{FKL} applications of wqh algebras to construct a unitary structure on a semisimple tensor category when it is equivalent to a unitary tensor category, and little information is needed on the tensor equivalence. Conversely,   we have shown that specific unitary
structures   on the weak Hopf algebra $A_W({\mathfrak g}, q, \ell)$ are  useful   to establish a unitary ribbon   equivalence   between the two tensor categories, the quantum group fusion category and the category of modules of the corresponding vertex operator algebra, in a clear    way.
This construction can be composed with     the previous step and provide a direct and explicit unitarization to module braided tensor categories of affine vertex operator algebras via quantum groups.

The specific structure alluded to and given in \cite{FKL} is that of 
  {\it discrete   unitary coboundary weak quasi-Hopf algebra.} This notion plays a central role to unify   Hermitian structures of   quantum groups at roots of unity and affine vertex operator algebras at positive integer level.
 This  definition is modeled on the  ribbon structure of $U_q({\mathfrak g})$ with its $^*$-structure for $|q|=1$.  
 
 We define the operator algebraic quantum group 
 as a unitary coboundary weak quasi-Hopf 
algebra. This structure is given by a discrete $C^*$-algebra $A$  equipped with the following data.

\begin{defn}\label{Hermitian_ribbon_wqh} A {\it Unitary coboundary} weak quasi-Hopf   algebra $A$ is defined by the following data:
\begin{itemize}
\item[{\rm        a)}]    A   weak quasi-Hopf  algebra $A$ endowed with a $C^*$-algebra
involution 
  with an antipode $(S, \alpha, \beta)$ 
\item[{\rm        b)}]   a
ribbon structure $(R, v)$ for $A$ associated to $(S, \alpha, \beta)$  
 such that the ribbon element $v\in A$
  is   unitary,  
\item[{\rm        c)}]     a unitary central square root $w\in A$
  of $v$ such that $\varepsilon(w)=1$, 
$S(w)=w$, 
\item[{\rm        d)}] $\tilde{A}=A^{{\rm        op}}$ as quasitriangular weak quasi-bialgebras,
\item[{\rm        e)}] $\overline{R}=Rw^{-1}\otimes w^{-1} \Delta(w)$ is positive.
 
  \end{itemize}
  
\end{defn}
  
  The operator $\overline{R}$ is Drinfeld coboundary matrix and gives the $\Omega$-involution $\Omega=\overline{R}$.
  The triple $(S, \alpha, \beta)$ describes the antipode of a weak quasi-Hopf algebra, and the presence of invertible elements $\alpha$, $\beta$ 
  is necessary in the quasi-Hopf cases, it arises from twist deformations of the usual antipode of a Hopf algebra.
The weak Hopf algebras $A_W({\mathfrak g}, q, \ell)$ are unitary coboundary. In this case $\alpha$ and $\beta$ are the identity of the algebra,
and $S$ satisfies the properties of the  usual antipode for Hopf algebras and commutes with the involution --  known as the Kac-type property in the setting of compact quantum groups.
The Zhu algebras $A(V_{{\mathfrak g}_k})$ become unitary coboundary weak quasi-Hopf algebras with a Drinfeld twist from $A_W({\mathfrak g}, q, \ell)$,  by our version of Drinfeld-Kohno theorem

  We have shown that a discrete unitary ribbon wqh algebra with unitary ribbon element and satisfying the compatibility property of the two opposite and adjoint $R$-matrices   is automatically unitary coboundary,
  and conversely,
 as one only has to choose a square root of the ribbon element, see Sect. 27 in \cite{FKL}.

 The construction of  the weak  Hopf algebras $A_W({\mathfrak g}, q, \ell)$ associated to
Wenzl's functor $W: {\mathcal C}({\mathfrak g}, q, \ell)\to{\rm Vec}$ for suitable primitive roots of unity $q$ is then a main result in our approach. When $q$ is a minimal root, we show that 
$A_W({\mathfrak g}, q, \ell)$  is 
unitary coboundary and its antipode is of  Kac type in a   sense motivated by the theory of compact quantum groups of Woronowicz.
This extends a result previously shown in \cite{CP} for the case ${\mathfrak g}={\mathfrak sl}_N$ with different methods.
In both cases, Wenzl's ideas on the fusion tensor structure generated  by the building tensor products $V_\lambda\otimes V$ are  an important starting point. Roughly speaking,   the decomposition into irreducibles is multiplicity free and allows to distinguish clearly the negligible representations, that belong to the tensor ideal defining the quantum group fusion category.   In the case of the vertex operator algebra, the corresponding representations, when lifted from the Zhu algebra to the whole vertex operator algebra, do not appear. On this basis, the  equivalence can be started and established by extending to all objects, keeping track of the respective unitary structures via the analogue of a Drinfeld twist. 

In the type $A$ case,  
one of our first results concerning Verlinde fusion categories in \cite{FKL}
 is Sect. 24 in \cite{FKL}, where we give a complete classification based on Kazhdan-Wenzl classification theory \cite{KW}. In this case,
 extending an argument by Neshveyev and Yamashita for the Hopf algebra $U_q({\mathfrak sl}_N)$ for generic parameter $q$, 
 to our weak Hopf algebra $A_W({\mathfrak sl}_N, q, \ell)$ for $q$ a root of unity, we show in particular that unitarity and braiding
 determine the associator completely.

Let us discuss extension to the other Lie types   in slight more detail. In these cases we use a direct method.
Let $q$ be a suitable complex root of unity.  Then the quantum group $U_q({\mathfrak g})$
in the sense of Drinfeld, Jimbo and Lusztig  has a $^*$-involution which is not a Hopf-algebra involution in the usual sense,
but  with the usual involution for real parameters it shares the property that the classical limit  reduces to the compact real form of ${\mathfrak g}$.
The coproduct satisfies an anticomultiplicativity relation
\begin{equation}\label{a_r}\Delta(a^*)=\Delta^{\rm op}(a)^*.\end{equation}
The quantum group has a highly non-trivial structure, the $R$-matrix  that was discovered by Drinfeld and Jimbo.
  This matrix gives rise to representations of the braid group in the representation category that with the involution
  satisfies the compatibility relation \begin{equation}\label{r-matrix_and_inv} R^*=(R_{21})^{-1}.\end{equation}
 We may interprete (\ref{a_r}) and (\ref{r-matrix_and_inv})   on the dual coordinate function algebra adopting Gelfand transform viewpoint.
Property (\ref{a_r}) corresponds 
  to the noncommutativity relation between operators
$(BA)^*=A^*B^*$
 and (\ref{r-matrix_and_inv})  corresponds to require that the induced braid group representations  on the adjoint  $(A^*B^*)^*$ and the opposite $BA$ multiplication structures coincide.
 In this way,  tensor product of unitary representations of   $U_q({\mathfrak g})$ is not $^*$-preserving with respect to the tensor product Hermitian form.
Similarly, the braid group representation on a tensor product of Hilbert space representations is not-unitary. To correct these features,
Wenzl showed in \cite{Wenzl} that
one has to change the inner product on tensor products.    A natural solution   arises from the ribbon structure itself. 
Indeed,   $U_q({\mathfrak g})$ is a ribbon Hopf algebra, in that the  square of the $R$-matrix is almost
the identity, and more precisely is given by a $2$-coboundary of a unitary central element $v$, also called the
quantum Casimir operator. The matrix $\overline{R}=R\Delta(v)^{1/2}v^{-1/2}\otimes v^{-1/2}$ can be given a rigorous
meaning at roots of unity and induces a non-degenerate Hermitian form with the desired properties, 
see   Sect. 30 in \cite{FKL}.
\smallskip

By the work of Andersen \cite{Andersen}, one has a notion of tilting module of $U_q({\mathfrak g})$, and a notion of negligible and non negligible module.
It is not easy to find a canonical choice of a fusion non-negligible submodule of a tensor product of non-negligible modules, see Sect. 11.3 in \cite{Chari_Pressley}. The negligible modules form an ideal in the category of tilting modules, and the quotient gives
an important   finite semisimple ribbon  category ${\mathcal C}({\mathfrak g}, q, \ell)$  with $\ell$ the order of $q^2$,
associated to $U_q({\mathfrak g})$,  which is modular and has a natural $^*$-structure introduced by Kirillov
 for some values of $\ell$. We refer to the review paper by Rowell for these results 
 \cite{Rowell2}, we also do some review work in Sects  19, 20, 30 in \cite{FKL}.
 \smallskip

  By a first  main result in \cite{Wenzl}, see also
Sect. 20 in \cite{FKL}, the Weyl modules $V_\lambda(q)$ of $U_q({\mathfrak g})$
have a unique invariant  positive inner product for $\lambda$   in the associated closed Weyl alcove. 
In particular, this holds for the objects in the open Weyl alcove, and these 
  correspond to the simple objects of ${\mathcal C}({\mathfrak g}, q, \ell)$.
In this way   the quotient category associated to the quantum group becomes a linear $C^*$-category.

A   main result of the same paper is the construction of a unitary tensor structure that involves at
the same time the construction of a canonical fusion tensor submodule $V_\lambda(q)\boxtimes V_\mu(q)$
of $V_\lambda(q)\otimes V_\mu(q)$.
 A central idea for the construction of a fusion submodule in \cite{Wenzl}   is the use of a generating representation $V$ such that  $V_\lambda(q)\otimes V(q)$ is completely reducible for $\lambda$ a simple object of
 ${\mathcal C}({\mathfrak g}, q, \ell)$  (except for the $E_8$ case  that is treated separately).

To make ${\mathcal C}({\mathfrak g}, q, \ell)$ into a unitary tensor category with unitary braided symmetry, one has to change the inner product on tensor products with a twist associated to the $R$-matrix and  a square root of the action of the ribbon element $v$, the quantum Casimir, on the tensor product space. 
It is shown in \cite{Wenzl} that  the Hermitian form    is positive definite on
Wenzl fusion submodules.  We review these construction in Sect. 30 of \cite{FKL}.\smallskip

This construction
  extends to all pairs of tensor products of simple objects and  
gives a natural $^*$-functor $$W: {\mathcal C}({\mathfrak g}, q, \ell)\to{\rm Hilb}$$ and natural
transformations $$F_{\lambda, \mu}: W(V_\lambda(q)\otimes W(V_\mu(q))\to W(V_\lambda(q)\boxtimes V_\mu(q))$$ 
$$G_{\lambda, \mu}: W(V_\lambda(q)\boxtimes V_\mu(q))\to W(V_\lambda(q))\otimes W(V_\mu(q))$$ such that $FG=1$.
The extension is implicit in \cite{Wenzl} and we make them explicit  using left paranthesized choices
of fusion subspaces
$$((V(q)\boxtimes V(q))\boxtimes V(q))\dots)$$
However, $G$ is not the adjoint of $F$ by non-triviality of the $R$-matrix, and in this sense we regard this structure
as non-unitary.\smallskip

  Wenzl showed using  Kashiwara bases that when $q_t$ varies continuously on the arc of the circle connecting $q$ to $1$ clockwise, then there are idempotents $p_\gamma(q_t)$ that in a suitable sense
  connect continuously the unitary representation $V_\lambda(q)\boxtimes V(q)$ of $U_q({\mathfrak g})$
  with a unitary representation $V_\lambda\boxtimes V$ of the classical Lie algebra with respect to the compact real form, with the same fusion
  rules and integral dimensions.\smallskip

 In the setting of vertex operator algebras, Zhu has associated an associative algebra $A(V)$ to a VOA $V$ that plays an important role
in the representation theory of  $V$. Under suitable rationality conditions due to Huang and Lepowsky that we recall in Sect. 21 of \cite{FKL}, the VOA $V$ has a semisimple representation category,  the Zhu algebra $A(V)$ is a semisimple algebra and $V$ and $A(V)$ have    equivalent representation categories as linear categories. The Zhu algebra may be defined as the algebra of natural transformation associated to the functor that takes a $V$-module $M$ to the top level subspace $M_0$.
It has long been asked how to construct a group-like structure on the Zhu algebra inducing a
tensor equivalence for their representation categories, see  the foundational paper by Frenkel and Zhu \cite{Frenkel_Zhu}. In  \cite{FKL} in the case of affine VOA at positive integer level we give such a 
natural construction, making $A(V)$ into a weak quasi-Hopf algebra. Indeed, by Tannakian duality the problem amounts
  to construct a weak quasi-tensor structure on Zhu's functor, and we construct this structure working with our natural weak tensor structure of the fibre functor associated to  $A_W({\mathfrak g}, q, \ell)$.
  
Another main result of \cite{FKL} is the formulation of an abstract Drinfeld-Kohno theorem showing that the class of unitary coboundary weak quasi-Hopf algebras
  satisfying an extra compatibility assumption between coproduct and involution is closed under the operation of
  taking a suitable square root construction of the associated positive operator $\Omega$ defining the unitary structure.
  Then we apply this result to the Drinfeld twist $\overline{R}\Delta(I)$ and we follow Wenzl continuous path and we have a weak quasi-tensor structure on Zhu's functor, but we shall come back to this twist.
Thus we have a direct connection from $A_W({\mathfrak g}, q, \ell)$ to the Zhu algebra $A(V_{{\mathfrak g}_{k}})$ of the affine Lie algebra $V_{{\mathfrak g}_{k}}$ for a suitable positive integer $k$.

When we assume that ${\rm Rep}(V_{{\mathfrak g}_k})$ is endowed with Huang-Lepowsky tensor structure then 
there is a tensor product module $M\boxtimes N$ for every pair of $V_{{\mathfrak g}_k}$-modules $M$, $N$.
By definition we have a bilinear map to a completed module
$$F: M\otimes N\to\overline{M\boxtimes N}$$ satisfying a universal property.
The associator of ${\rm Rep}(V_{{\mathfrak g}_k})$
is related to  KZ equations by their work. 

Composing with inclusion and projection onto the top level modules, Huang-Lepowsky theory gives 
bilinear maps
$$F^{HL}_0: M_0\otimes N_0\to (\overline{M\boxtimes N})_0.$$
We endow ${\rm Rep}(V_{{\mathfrak g}_k})$ with the Hermitian structure that can be regarded in two ways, either directly from
Kirillov work or as the form induced by Wenzl continuous path and our twist. 
We show  that   $F^{HL}_0$ can be completed to a weak quasi-tensor structure $(F^{HL}_0, G^{HL}_0)$  in a canonical way, that is 
$$F^{HL}_0G^{HL}_0=1.$$
In particular, on the generating pairs  $(V_\lambda, V)$, with $V$ the fundamental representation of the Lie algebra
and $V_\lambda$ any simple representation, the structure satisfies the unitarity condition
$$G^{HL}_0=(F^{HL}_0)^*, \quad\quad F^{HL}_0(F^{HL}_0)^*=1.$$

We review this tensor product for models
that may be compared to the WZW model, in an abstract form, and focusing on the structure of the lowest energy subspace.  

  In this way the tensor product structure of Huang and Lepowsky and Tannaka-Krein duality
  applied to Zhu's functor with the weak quasi tensor structure $(F_0, G_0)$ make   $A({V_{\mathfrak g}}_k)$
  into a wqh algebra in a natural way.
  In other words,  $A({V_{\mathfrak g}}_k)$ becomes a   weak quasi-Hopf algebra with representation category ribbon equivalent to
${\rm Rep}(V_{{\mathfrak g}_k})$.
Moreover, 
  $A({V_{\mathfrak g}}_k)$ is a
 unitary coboundary wqh algebra isomorphic up to twist and isomorphism to
 $A_W({\mathfrak g}, q, \ell)$, and Kirillov Hermitian form
 makes ${\rm Rep}(V_{{\mathfrak g}_k})$
into a unitary modular tensor category derived directly from a structure with the same property from the quantum group fusion category, see Sect. 7, 33 of \cite{FKL}.
We  study   connections with the braiding structure
arising from loop group fusion categories in the work of A. Wassermann, Toledano-Laredo, Gui in Sect. 34 in \cite{FKL}.

Our construction makes
$A(V_{{\mathfrak g}_k})$ into a unitary coboundary weak quasi-Hopf algebra in a natural way with the structure
transported from the quantum group. In this way  we give
   a   direct proof of Finkelberg-Kazhdan-Lusztig theorem in Theorem 2.4 in \cite{FKL} in the setting of vertex operator algebras with Huang-Lepowsky braided tensor structure. This is done in  Sect. 33--38 of \cite{FKL} based on the construction of $A_W({\mathfrak g}, k, \ell)$ of the previous section, and taking into account the last part of proof given in \cite{Pinzari_constructing} that is based on the current literature on quantum Schur-Weyl duality for
   the choice of a fundamental representation for each Lie type.

  Along the way, 
  we show that by general reasons due
to amenability properties, a unitary weak Hopf algebra   with trivial unitary structure (coproduct and involution commute)
and amenable representation category has necessarily integer dimensions. This extends previously known facts for Hopf $C^*$-algebras.  
It follows that    the untwisted algebra has the following alternative.
Either the involution commutes with the coproduct, and then it automatically has an "exponential" $R$-matrix and unitary structure, and then it must have a $3$-coboundary non-trivial associator (Drinfeld-Kohno model)
or else one can not obtain the commutation relation between coproduct and involution globally, but only on some
possibly sufficiently many pairs of representations, and this motivated part b) of our Drinfeld-Kohno theorem 29.4 in \cite{FKL}. In this case, the weak  Hopf structure, that is triviality of the associator in the weak sense, is not ruled out.
 In our case the alternative applies to the Zhu algebra, for which we have indeed constructed  
 the structure of a unitary coboundary weak-quasi-Hopf algebra  with a $3$-coboundary associator
associated to the corresponding affine VOA following Wenzl continuous path recalled in the abstract, and inducing   Huang-Lepowsky tensor structure, see also
Sect. 25--34 in \cite{FKL}.
We show that on a sufficiently large collection ${\mathcal V}$ of representations, this associator has the form of the associator
of a weak  Hopf algebra. We discuss general aspects on the associator in Sect. 8 and complete the proof in \cite{Pinzari_constructing}.

  We come back to some more detail on the unitary structure of ${\mathcal C}({\mathfrak g}, q, \ell)$ of \cite{Wenzl} that reflect on $A_W({\mathfrak g}, q, \ell)$.
Wenzl considers two main problems. Since not all the specialized Weyl modules $V_\lambda(q)$ are irreducible and complete reducibility does not hold, we have that    not all
the idempotents $p_\gamma:   V_{\lambda, {\mathcal A}'}\otimes_{{\mathcal A}'} V_{\mu, {\mathcal A}'}\to V_{\gamma, {\mathcal A}'}$ describing classical fusion
    specialize to corresponding idempotents $p_\gamma(q)$. This is an obstruction in 
  the construction of positive Hermitian structures on tensor products of two specialized
Weyl modules $V_\lambda(q)\otimes V_\mu(q)$  associated to the braiding as before.

Indeed, when  the obstruction vanishes, then the Hermitian form of the specialized
    module $V_\lambda(q_0)\otimes V_\mu(q_0)$ at a fixed suitable root of unity $q_0$ associated to the braiding is positive on $V_\gamma(q_0)$ by a continuity argument that
    links $q_0$ to the classical limit $1$ that  comes from specialization.
    Notably, the building tensor products $V_\lambda(q_0)\otimes V(q_0)$ turn out to be   sufficiently many   for which this holds. From these,    the weak  Hopf $A_W({\mathfrak g}, q, \ell)$ can be constructed,  and as already mentioned is unitary coboundary when $q$ is minimal Sect. 29 in \cite{FKL}.

However it is not clear whether  the ribbon element (the quantum Casimir) $v$ admits a square root
  in a reasonable  topological completion
of the integral form ${\mathcal U}_{{\mathcal A}'}^\dag({\mathfrak g})$. This problem is mentioned
at the beginning of Sect. 3.6 in   \cite{Wenzl} and is related to the construction of an invariant Hermitian form of the full tensor product spaces $V_\lambda(q_0)\otimes V_{\mu}(q_0)$. More precisely, note that if   such a square root   lied in a   topological completion then
by our   theorem 27.13 in \cite{FKL} we would have an Hermitian form on the full tensor
product obtained by specialization of a nondegenerate  ${\mathcal A}'$-valued Hermitian form.
Since it is obtained as a specialization, following Wenzl, the Hermitian form  would be continuous in $q$, and thus positive since it reduces to the
usual inner product for $q=1$. It would then follow that the tensor product is completely reducible, but this is not always the case.
 This problem
has been overcome in \cite{Wenzl} for
  tensor products that decompose into a direct sum of specialized Weyl modules
  (except for the $E_8$-case, where complete reducibility does not hold, but he solves the problem in this case in a direct way).
  Then he uses some of them which arise as tensor products of a fundamental representation $V$ of the Lie algebra with any other
  irreducible representation $V_\lambda$ in the open Weyl alcove.

For our purposes,  we   need    an Hermitian  form associated to Drinfeld coboundary on all tensor products $V_\lambda(q)\otimes V_\mu(q)$ of pairs of irreducibles
 in the open alcove. Relying on the Tannakian approach   of Sawin to $U_q({\mathfrak g})$, we shall propose a 
solution that regards the action of $v$  on a tensor product as an operator. 
We  take a canonical square root $\Delta(v)^{1/2}$
non-formally as an operator on a tensor product of two specialized Weyl modules in the open alcove. 
Since $R_{21}R$ is the
$2$-coboundary of $v$, 
this way we have a square root of  $R_{21}R$
and we
  construct an action of
Drinfeld coboundary matrix
  $\overline{R}=R(R_{21}R)^{-1/2}$ on the full tensor product, which turns out selfadjoint and invertible.
  Following \cite{Wenzl}, with respect to the associated Hermitian form,  tensor products of Hermitian (e.g. Hilbert space)
   $^*$-representations of $U_q({\mathfrak g})$ is Hermitian and
  the braided symmetry is unitary.

These arguments are in Sect. 19, 20, 30 in \cite{FKL}.
After doing this, we have a nondegenerated Hermitian form defined by the action of $\overline{R}$,
and we consider the twist $F=\overline{R}^{1/2}\Delta(I)$ supported on a domain that arises from a specialization,
and thus continuity arguments may be applied, and lead to an application    of our abstract Drinfeld-Kohno theorem
   in   Sections 20, 31, 33, 34.

  \section{Directions of future research}

 We conclude with some   research lines concerning the specific unitary structures of wqh.    \smallskip

\noindent{\bf Problem of reconstruction of the field algebra from the observable algebras for the WZW model.}   Mack-Schomerus original motivation for the introduction of weak quasi-Hopf algebras
was the construction of the field algebra from the local observable algebra in the algebraic approach
to CFT by Doplicher-Haag-Roberts \cite{MS, MS1, Schomerus_field_algebra}  generalizing Doplicher-Roberts corresponding construction for compact groups in high dimensional AQFT \cite{DHR}.
The field algebra has been reconstructed by Schomerus in \cite{Schomerus_field_algebra}. In this paper the weak quasi-Hopf algebra was constructed using only a weak dimension function, and this leads to highly non-unique weak quasi-Hopf algebras, which were not canonically associated to the category arising from the conformal field theory. Thus a main problem remains that of constructing canonical quantum gauge group to which apply Schomerus construction of the field algebra, to study.

In \cite{FKL} we have given   a categorical characterization of unitary coboundary weak quasi-Hopf algebras.
This characterization reduces for braided symmetric tensor categories to the notion of symmetric tensor functor,
which plays a central role in Doplicher-Roberts duality theory for compact groups \cite{DR_CMP}, \cite{DR1}. It is natural to ask whether
in the case where we have a unitary coboundary weak  Hopf algebras or a unitary coboundary wqh algebras
with a $3$-coboundary associator,  e.g. arising from ${\mathcal C}({\mathfrak g}, q, \ell)$ or affine VOAs at integer level, these constructions may be further developed to study Schomerus field algebra in detail as induced from an ergodic action of a weak Hopf algebra
extending our previous results summarized in the introduction in the setting of compact quantum groups. \smallskip

\noindent{\bf Problem on Reshetikin-Turaev invariants}  
  Reshetikhin and Turaev    introduced
  ribbon and modular Hopf algebras   to describe Jones and Witten invariants of knots, links and $3$-manifolds. In their work, they needed to use non-semisimple
  Hopf algebras associated to  semisimple categories arising from conformal field theory   \cite{RT}, \cite{Turaev}.

We have developed  an abstract Drinfeld-Kohno theorem that gives a canonical unitary ribbon equivalence between unitary coboundary wqh algebras and thus preserves modularity properties. In particular, the categories have the same link and manifold invariants.
We prove that discrete Hermitian (unitary) coboundary wqh are   Hermitian (unitary) ribbon   if and only if the antipode satisfies a certain condition, that is automatic when the antipode is of Kac type.
 Thus a discrete unitary coboundary wqh   leads to Reshetikhin-Turaev  link invariants  \cite{Turaev}.

 We have developed a general theory of Hermitian coboundary wqh algebras.
 For the application we have mostly focused on the unitary cases corresponding to the "minimal" roots
 of unity $q=e^{i\pi/\ell}$ with $d|\ell$, where $d$ is the ratio of the square of the lengths of the long to the short roots. These values  correspond
 to   the cases of most physical interest  of positive integral levels.
 However, most of our results, that is construction of the weak  Hopf algebra hold for non necessarily minimal   roots of unity.
 In these cases,  we loose unitarity of the structure, but we still have Hermitian structures.
 We show that Hermitian (discrete) coboundary wqh algebras are always Hermitian ribbon wqh algebras in Sect. 27 in \cite{FKL}.
 Moreover, these cases might still be covered by our Drinfeld-Kohno theorem 29.4 in the weaker form
 (that is dropping the assumptions in b) and c) in the Hermitian case). The corresponding categories may overlap with those considered by Rowell
\cite{Rowell1}, \cite{Rowell3}. Thus they possibly correspond to Chern-Simons theory with fractional level,
see the last remark in \cite{Sawin}. It is natural to ask whether the notion
of  discrete Hermitian coboundary wqh with compatible involution may be useful to the study
  modular Hermitian of unitary modular categories, TQFT, link and manifold invariants.
  
  It is also natural to ask whether there is a relation between the discrete coboundary property and modularity,
 or whether discrete coboundary wqh algebras lead to  invariants of closed oriented $3$-manifolds, or whether there are connections with
 Witten invariants. 
 We propose to study the connections between Reshetikhin-Turaev and Witten invariants,
   as   invariants of  our weak (quasi)-Hopf algebras regarded as quantum quotients of $G$ directly associated to the conformal field theory.
   \smallskip
 
\noindent{\bf   Problem on the search of Dirac operators} A further   direction of future research   arises from noncommutative geometry in the sense of Connes
\cite{Connes}. Focusing attention to the WZW model of conformal field theory associated to a compact Lie group
$G$, the Dirac operator on $G$ may give rise to a Dirac operator on the "quantum quotient" of $G$ describing the conformal field theory,
that is a very interesting  subject to study.

One of the main results of  \cite{FKL} is the construction of a quantum quotient of $G$
for WZW model at positive integer level as a canonical wqh structure on the Zhu algebra
$A(V_{{\mathfrak g}_k})$  associated to Huang-Lepowsky tensor product theory and Kirillov-Wenzl unitary structure.
 
Neshveyev and Tuset    used a Drinfeld twist   to construct a Dirac operator on
the compact quantum group associated to $U_q({\mathfrak g})$ 
connecting this quantum group to the Drinfeld category associated to $U({\mathfrak g})$, for $0<q<1$ 
\cite{NT_twist}. 

Given that we are following their scheme
when $U_q({\mathfrak g})$ is replaced by our weak Hopf algebra $A_W({\mathfrak g}, \ell, q)$,   $U({\mathfrak g})$ is replaced by the Zhu algebra,
the representation category of $U_q({\mathfrak g})$ is replaced by the quantum group fusion category, the Drinfeld category is replaced by the affine VOA module category,
it is natural to ask to study   the search of Dirac operators associated to a Casimir operator on
our weak  Hopf algebras $A_W({\mathfrak g}, \ell, q)$ and, following Neshveyev-Tuset method of Drinfeld's twist.
on the Zhu algebra 
of   the affine Lie VOA $V_{{\mathfrak g}_k}$.   

Moreover, Dirac operators in Conformal Field Theory have been also studied
with different methods by other authors. It seems a natural and interesting problem to study connections between the various approaches.

\medskip

\noindent{\bf Problem of extending our correspondence beyond the affine vertex operator algebras. Study of extensions of VOAs and conformal nets.} One more motivation may arise from the problem of comparing different  axiomatizations of conformal field theories. A geometrical approach has been initiated by Segal \cite{Segal}. Then as already mentioned we have 
an  algebraic approach based on vertex operator algebras \cite{Borcherds}, \cite{FLM} and a functional analytic approach in the setting of conformal nets \cite{Haag}. A general connection 
has been studied between the vertex operator algebra and the conformal net approach under certain circumstances, and one would like to
extend these results to more general situations, and extend the theory to representation categories,
reconciling with
the field algebra construction of Doplicher and Roberts \cite{DR_CMP} in high dimensional AQFT
 from a   common quantum gauge group.
 We have recalled in Sect. 1 that in the case of rigid tensor $C^*$-categories with no braided symmetry, the problem has been studied e.g. 
 in \cite{PRergodic, PRinduction} in the setting of compact quantum groups. A next step is the study of the Virasoro algebra from the view point of quantum homogeneous spaces.
Another natural problem is that of a comparison of our methods  with the study of extensions of conformal nets and vertex operator algebras in the literature.
 \medskip

\noindent{\bf Problem on theory and construction of quantum quotient spaces for weak Hopf algebras.} This problem is closely related to the previous problem. It is   natural to ask whether one can develop  a theory of quotient spaces in the sense of noncommutative geometry for wh (or even wqh) algebras in some analogy with the case of compact quantum groups, and construct examples associated to fusion categories of CFT which do not admit a natural weak
 dimension function. Natural examples would be the Virasoro VOAs,
   some more examples are discussed in Sect. 22 in \cite{FKL}. We recall that in the case of compact quantum groups,
 quotient spaces enjoy a natural Tannakian duality, and this has been the subject of   investigation
 in connection to noncommutative geometry   see e.g. \cite{Podles_quantum_spheres}, \cite{PRergodic}, \cite{DeCommer_Yamashita},   \cite{DNTY1}, \cite{DNTY2}.\smallskip

\noindent{\bf Problem of comparison with Gui's hermitian form.} Gui has defined an invariant Hermitian form and showed positivity for a large class of VOAs, including the affine cases
extending work of Wassermann.
 It is natural to ask whether  Gui  Hermitian structures coincides with our unitary structure  on the Zhu algebra.\smallskip

\noindent{\bf Beyond $\textbf{2D}$.} As suggested in the introduction, the braided Doplicher-Roberts program faces significant challenges when situated within the broader context of modern QFT.

 Recent  research by Mund, Rehren, and Schroer shows the emergence of the braid group even for $4D$ in QED, 
for the description of charged particles
that are string-localized \cite{MRS}.
See also \cite{Wuang_Hazzard}, and \cite{Zhou} and references therein for recent results and discussion on parastatistics. 

Furthermore, the connection between   higher-dimensional   and lower-dimensional theories has been deepened by Rehren's work on {\it algebraic holography}  \cite{Rehren00} in the setting of AQFT and by the {\it celestial holography program}  in physics
where the asymptotic symmetries of $4D$ Yang-Mills theory are related to $2D$ conformal symmetries  governed by the WZW model  \cite{HeMitraStrominger16}.

 These developments are particularly significant given that the rigorous construction of the local observable algebras and field operators for $4D$ Yang-Mills theory remains one of the most profound open problems in mathematical physics.  While the original $4D$
 Doplicher-Roberts theory provides a rigorous reconstruction of the gauge symmetries and the field algebra under the setting of permutation statistics, the extension of this framework to   braided categories   found in low dimensional QFT ($D\leq3$),
or in these modern holographic contexts, or the braid group exchange relations satisfied by string-localized fields in the algebraic approach to QED of \cite{MRS} remains a fundamental challenge for the braided Doplicher-Roberts program.

 \end{document}